\documentclass[1p,10pt,a4paper]{elsarticle}

\usepackage{lineno,hyperref}
\modulolinenumbers[5]
\usepackage{latexsym}
\journal{Journal of \LaTeX\ Templates}
\input{amssym}

\usepackage{numcompress}

\newtheorem{thm}{Theorem}
\newtheorem{lem}[thm]{Lemma}
\newtheorem{dfn}{Definition}

\newtheorem{example}{Example}
\newtheorem{cor}[thm]{Corollary}
\newdefinition{rmk}{Remark}
\newproof{pf}{Proof}
\begin{document}

\begin{frontmatter}

\title{Hartman-Grobman Theorem for IFS}

\author{Mehdi Fatehi Nia\corref{fn1}}
\cortext[fn1]{Corresponding author}
\ead{fatehiniam@yazd.ac.ir}
\address{PhD, Department of Mathematics, Yazd University, Yazd 89195-741, Iran}

\author{Fatemeh Rezaei\corref{first}}
\cortext[first]{Principal corresponding author}
\ead{f$\_$rezaaei@yahoo.com}

\begin{abstract}
In this paper, for iterated function systems, we define the classic concept of the dynamical systems: topological conjugacy of diffeomorphisms. We generalize the Hartman-Grobman theorem for one dimensional iterated function systems on $\Bbb R$. Also, we introduce the basic concept of  structural stability for an iterated function system and so we investigate the necessary condition for structural stability of an iterated function system on $\Bbb R$.
\end{abstract}


\end{frontmatter}

\section{Introduction}
This section includes three parts. In the first part, there is provided an almost perfect review of the literature on studies which have been done on the iterated functions systems. Also we introduce their applications to understand the importance of studying of the IFSs. In the second part, we describe the history of the creation and importance of one of essential theorems of the local dynamic that is named Hartman-Grobman Theorem and in the following, we study the researches done on the generalization and extension of this theorem. In the third part of this section, we define a very important concept of dynamic, which is related to this theorem, and is named structural stability. Moreover, we briefly state the history of presentation of this concept.

\paragraph{First part}
The concept of the iterated functions systems was applied in 1981 by Hutchinson. Moreover, the mathematical basic of the iterated functions systems was established by him;\cite{HUTCHINSON1981}, but this term was presented by Barnsley; briefly as IFS;\cite{HS2012}.
We know that an IFS includes a set $\Lambda$ and some functions $f_\lambda,\lambda\in \Lambda$, on an arbitrary space $M$. As, in an IFS, the set $\Lambda$ can be finite or infinite(countable) or its functions can be special, so different IFSs have been investigated. The most studies on the finite IFSs have been done by Barnsley;\citep{BARNSLEY1985,BARNSLEY1988,BARNSLEY1993,BARNSLEY2006,ABVW2010}. We can see the generated countably IFSs in some articles like\cite{MAULDIN1996}. In \cite{ABVW2010}, the IFSs have been studied, whose functions are affine transformations on the Euclidean spaces; there can be found many studies on this case of the IFSs. Also, in \cite{BV2012} the IFSs whose functions are onto transformations on the real projective spaces were investigated. Until recently, this subject has extended to the functions that are Mobius transformations on complex plane( or equivalently Riemann sphere);\cite{VINCE2013}. Although generally in \cite{Younis2011} the IFSs on manifold have been studied, the IFSs can be found that are considered on Hilbert spaces, and are called Perry IFSs;\cite{PERRY1996}. Thus, the IFSs were extended to conformal IFSs;\cite{mauldin1999}. We can see the generalization of Barnsley's concept at \cite{OlSEN2007},\cite{OLSEN2008} and \cite{OLSEN20082}.
We know that attractor of an IFS is something that is said fractal. But what is a fractal? We have understood that we can not have accurate description of geometric structure of many natural things like clouds, forests, mountains, flowers, galaxies and so on by using classical geometry. Mandelbrot, 1982, changed this perspective through which classical geometry extended into, so called, fractal geometry. In fact, fractal is made of iterating the functions in a set that is called iterated functions system. The literature have been proved that a fractal is the attractor of an IFS. The IFS model is a base for different applications, such as computer graphics, image compression, learning automata, neural nets and statistical physics\cite{EDALAT1996}. So, the study of the fractal is important and therefore, from one point of view, the study of an IFS as the way that can generate a fractal, is important;\cite{BARNSLEY1985}. The existence and uniqueness of the attractor of a finite IFS in 1985 was proved by Hata in \cite{hata1985}, also you can see \cite{duvall1992}. Abundant studies have been fulfilled on the context of the topological properties(such as dimension, measure, separation property) of an attractor;for example \citep{MAULDIN1996,chin1997correlation,mauldin1999,Nguyen2002,dumitru2011topological,barnsley2014}. The attractor of the affine IFS has many applications; for example image compression\citep{Jacquin1992,Reusens1994,FISHER1995,barnsley1996,Barnsley1988image,Zheng2006improved}, geometric modeling \citep{Blanc1997,stotz2008ifs,Alexander2012}. Moreover, the IFSs that are said recurrent IFSs \cite{barnsley1989recurrent}, have applications at generation of digital images because these images have curves that are not generable using standard techniques;\cite{barnsley1988application}. Also, the IFSs can be used as tools for filtering and transforming digital images;\cite{barnsley2011transform}. From a dynamic point of view, the IFSs have been studied, for example the stability and the hyperbolicity of an IFS in \cite{Ambroladze1999stability} and the asymptotic stability of a countable IFS( presentation of sufficient condition) in \cite{Janoska1995stability}, the asymptotic behaviour of a finite IFS with contraction and positive, continuous, place dependent probabilities functions in \cite{Jaroszewska2002}. Now, we also study a dynamic property, structural stability, for the finite IFSs that have not been investigated yet.

\paragraph{Second part}
Now we start writing the second part by presenting the concept of "topologically conjugate".
Sometimes, it has been seen that two systems seemingly are different but if we investigate these systems, dynamically, we find that they have the same behavior;\cite{Meiss2008}. In other words, the two systems are "equivalent", that is, studying one system will provide dynamic information about the other system. Thus, such systems allow us to look for some(approximately) simple system or an identified one equivalent to the complicated system to study that. In this context, concepts were proposed, which called topological equivalence and topological conjugacy. In this paper, we also define these concepts for IFSs and determine a special class of IFSs that has such properties.\\
A fundamental theorem to study the local behavior of a system that is a strong tool in dynamical systems is well-known as the Hartman-Grobman theorem or linearization theorem. This theorem examines the local behavior of a system about hyperbolic fixed points; accurately, the theorem say that the dynamical behavior of a system is the same as the dynamical behavior of its linearization near the hyperbolic fixed points. Thus, we can locally draw the phase space about these special points. Specially, this matter is important when the given system is nonlinear;\cite{Zimmerman2008}. Formation of this critical theorem is a question asked by M. M. Peixoto as following:\\
Consider the real, autonomous and nonlinear system
$Y^{'}=AY+E(Y)$
where $Y$ is a vector, and $\mid Y\mid$ is the Euclidean length of $Y$ and $E(Y)$ a smooth vector valued function for small $\mid Y\mid$ and also $E(Y)=o\big(\mid Y\mid\big)$ as $Y\rightarrow 0$ and  matrix $A$ is constant with the eigenvalues whose real parts are not zero. The question is whether there exists a topological mapping like $h(u)$ that $u=u(Y)$ from a neighborhood $Y=0$ onto a neighborhood $u=0$ such that solution paths of the above system locate in trajectories of solution $u^{'}=Au$?; \cite{Hartman1960lemma}\\
In 1959, Hartman answers the question in \cite{Hartman1960lemma}. In fact, he proved that if $E$ belongs to $C^2$-class then the answer of this question is positive. Also, in \cite{Hartman1960}, he has demonstrated that even if the function $E$ is analytic, the function $h$ does not need to be of class $C^1$. Also, according to the literature, including \cite{Hasselblatt}, Grobman, separately (maybe) in \cite{Grobman1959} has provided demonstration in 1959, therefore this theorem is well-known as the Hartman-Grobman Theorem, (briefly as H-G-T). Moreover, in 1963, Hartman proved that if, for all small $Y$, the function $E$ belongs to the class $C^1$ (or $F$ is uniformly continuous Lipschitz such that the constant of Lipschitz approaches to zero as $Y\rightarrow 0$), then the function $h$ also exists; that is, the answer of the Peixoto's question is positive;\cite{HARTMAN1963}.\\
In 1968, J. Pails extended the H-G-T for maps to infinite Banach space and for this case, gave a short proof;\cite{Palis1968}. In 1969, Pugh using Moser's techniques in \cite{MOSER1969} for this case rendered proof in \cite{PUGH1969}.
But general state of the H-G-T for maps in Banach spaces was proved by Quandt in \cite{QUANDT1986}.
The H-G-T enlarged non-autonomous systems. Palmer in 1973 generalized this theorem for non-autonomous systems whose linearization have exponential dichotomy;\cite{PALMER1973}.
The Linearization also reached random dynamical systems(shortly RDS) and control systems. The first linearization of discrete and random dynamical systems is related to Wanner's work in 1994;\cite{wanner1995linearization}.
The H-G-T for discrete and random dynamical systems in 2003 was generalized by Coayla-Teran and Ruffino in \cite{Coayla2003}.
These authors and A. Mohammed extended this theorem for continuous and random dynamical systems in 2007; we can see the results in \cite{Coayla2007}. Moreover, in this article, the H-G-T was generalized for hyperbolic stationary trajectories at RDS. In fact, it was proved that conjugate relation exists between the trajectories in the neighborhood of the origin and the corresponding neighborhood at tangent space.\\
The linearization theorem was expanded for controlling systems with special inputs by Pomet, Chyba and Baratehart in 1999 and was proved in \cite{baratchart1999}.\\
The H-G-T was enlarged for non-autonomous systems with discrete time. As, the linearization theorem has applications, specially, at theme partial differential equations on Banach spaces, the extension of H-G-T on Banach spaces is important. We can name some researchers as L. Barreira  and C. Valls in 2005, who did researches on  Banach spaces for non-uniformly hyperbolic dynamic;\cite{Barreira2006} and they extended this theorem for very general non-uniformly  hyperbolic dynamic;\cite{Barreira2009}.
Sola-Morales and M. Rodrigues generalized the H-G-T for infinite spaces with special conditions such as Hilbert space;\citep{RODRI2004,RODRI}. Another research done on Banach spaces in order to expand the H-G-T for maps was by V. Rayskin and G. Belitskii in\cite{Belitskii2009};they showed that at special conditions the homeomorphism is $\alpha$-Holder.\\
We know that the Hartman-Grobman theorem says that any $C^{1}$-diffeomorphism is topologically equivalent to its linear part at neighborhood of the hyperbolic fixed points. We also generalize this theorem for IFSs. In fact, we show that if the origin is hyperbolic fixed point of the $C^{1}$-diffeomorphisms of IFSs $\mathcal{F}$ and $\mathcal{G}$ and all the derivatives of these functions at zero belong to the same interval $(0, 1)\Big($ or $(-1, 0)$ or $(1, +\infty)$ or $(-\infty, -1)\Big)$, then these two IFSs have the same dynamical behaviors.

\paragraph{Three part}
Though, sometimes, systems look like, seemingly, they have completely different dynamical behaviors,(it raises bifurcation, chaos,...). Therefore, it leads to creating another concept that is called "structural stability". The literature such as \cite{bonatti2001dynamical} said that the concept of the structural stability with this name was introduced by M. M. Peixoto. In fact, this concept is the generalization of the concept of systems grossier or rough systems in 1973 by A. A. Andronov and L. S. Pontryagin. Andronov was interested the preservation of the qualitative properties of the flows under small perturbations and asked a question whose history can be seen in \cite{Anosov1985}. Indeed Peixoto in 1959 introduced the concept of the structural stability using corrections of the mistakes of the article \cite{Baggis1955}.
 We say that $C^{k}$-diffeomorphism $f$ is structurally stable if there exists a neighborhood of $f$ in the $C^{k}$-topology such that $f$ is topologically conjugate to every function at this neighborhood; in accurate words, a $C^{k}$- diffeomorphism $f$ is structurally stable if  for any $\epsilon>0$ there is a neighborhood $U(\epsilon)$ of $f$ in the $C^{k}$-topology such that any $C^{k}$- diffeomorphism $f_1\in U(\epsilon)$ is topologically conjugate to $f$, \cite{Pilyugin1999}. We consider the distance between two IFSs as maximum distance between the functions of two IFSs and so nearby IFSs makes sense. Thenceforth, we define the concept of structural stability for IFSs. Moreover, we demonstrate that necessary condition for an IFS to be structural stability is that all the fixed points of IFS' functions should be hyperbolic.

\section{A Preliminary Lemma}
We know that the diffeomorphisms $f,g\,:\Bbb {R}^{m} \rightarrow \Bbb {R}^{m}$ are topologically conjugate if there exists a homeomorphism $h\,:\Bbb {R}^{m} \rightarrow \Bbb {R}^{m}$ such that $hof=goh$, or equivalently $f= h^{-1}ogoh$. The function $h$ be said topological conjugacy.
Also for given $\epsilon>0$, we say that the diffeomorphisms $f$ and $g$ are $\varepsilon$-topologically conjugate if there exists a topological conjugacy $h$ such that $\|x-h(x)\|<\epsilon$ for every $x\in \Bbb {R}^{m}$. ($\|.\|$ is the norm on $\Bbb {R}^{m}$) \cite{Pilyugin1992introduction}\\
The following lemma has been proved at the some literature like \cite{Palis1980}(it has been given at \cite{Devaney1989introduction} as a practice). But as this lemma has critical role in demonstration some of theorems of this paper, we give proof with complete details, also it is remarkable to see these details. This lemma say that the contractive functions on $\Bbb{R}$ are topologically conjugate.
\begin{lem}\label{functions}
Suppose that real value functions $f$ and $g$ on $\Bbb R$ are defined with criteria $f(x)=kx$ and $g(x)=mx$ where $0<k,m<1$. Then $f$ and $g$ are topologically conjugate.
\end{lem}
\begin{pf}
Let $a$ be arbitrary positive real number. We know that there exists a homeomorphism $h$ from $[f(a),a]$ to $[g(a),a]$ such that $h(f(a))=g(a)$ and $h(a)=a$. Suppose $x\in\Bbb R$ is arbitrary and more than $a$. As for criterion of function $f$, as $n$ increases, the value $f^{n}(x)$ approaches the origin, as a result, there exists $n\in\Bbb N$ such that $f^{n}(x)<a$, assume that $n_x$ is the first $n$ with this property; that is, $k^{n_x}x<a<k^{n_x-1}x$. By considering the inequality $a<k^{n_x-1}x$, we have $ka<k^{n_x}x$, so $ka<k^{n_x}x<a$; that is, $f^{n_x}(x)= k^{n_x}x\in[f(a), a]$. Now for every $x>a$, we define the function $h$ as follows\\
\begin{equation}
\left\{ \begin{array}{rl}
 h:(a,+\infty) & {\longrightarrow (a, +\infty)}\\
 x & \longmapsto {g^{-n_x}\Big(h\big(f^{n_x}(x)\big)\Big)}
\end{array}\right.
\end{equation}
Firstly, $h$ is well-defined by considering the way $n_x$ is chosen. Moreover, the range of the function $h$ is
$(a, +\infty)$ since for every $x>a$, $f^{n_x}(x)\in[f(a), a]$, hence, in the basis of the definition of the function $h$ and its continuity we have $h\big(f^{n_x}(x)\big) \in[g(a), a]$; that is, $ma=g(a) < h\big(f^{n_x}(x)\big) < a$, and since the function $g^{-n_x}$ is strictly increasing, we obtain the following relation:
\begin{eqnarray*}
g^{-n_x}(ma) < g^{-n_x}\Big(h\big(f^{n_x}(x)\big)\Big)
\Longrightarrow({1\over m})^{n_x}.ma < g^{-n_x}\Big(h\big(f^{n_x}(x)\big)\Big)\\
\Longrightarrow({1\over m})^{n_x-1}a < g^{-n_x}\Big(h\big(f^{n_x}(x)\big)\Big),
\end{eqnarray*}
We know $n_x\in\Bbb N$ and since $0<m<1$ so ${1\over m} > 1$, and consequently, we have $a<g^{-n_x}\Big(h\big(f^{n_x}(x)\big)\Big)$; that is, $h(x)>a$.
Secondly, the function $h$ is a homeomorphism because it is a composition of the homeomorphisms.\\
Now, suppose $0<x<f(a)$; that is, $0<x<ka$ so $0<{1\over k}x<a$ and this means that $0<f^{-1}(x)<a$ then based on the criterion of the function $f^{-1}$, as $n$ increases, the value of $f^{-n}(x)$ gets far away the origin, so there exists $n\in\Bbb N$ such that $f^{-n}(x)> f(a)$. Assume $n_x$ is the first $n$ with this property; that is,
$({1\over k})^{n-1}x < ka < ({1\over k})^{n}x$. By considering the inequality $({1\over k})^{n-1}x < ka$ we obtain $({1\over k})^{n}x < a$, so $ka < ({1\over k})^{n}x < a$; that is, $f^{-n}(x)= ({1\over k})^{n}x\in [f(a), a]$.\\
For every $x\in\Big(0, f(a)\Big)$ we define the function $h$ from $\Big(0, f(a)\Big)$ to $\Big(0, g(a)\Big)$ with the criterion $h(x)= g^{n_x}\Big(h\big(f^{-n_x}(x)\big)\Big)$. Clearly, this function is well-defined and the range of the function $h$ is $\Big(0, g(a)\Big)$ since for every $x\in\Big(0, f(a)\Big)$, $f^{-n_x}(x) \in [f(a), a]$ so
$h\big(f^{-n_x}(x)\big) \in [g(a), a]$; that is, $ma=g(a) < h\big(f^{-n_x}(x)\big) < a$, and because the function $g^{n_x}$ is strictly increasing then
\begin{eqnarray*}
g^{n_x}\Big(g(a)\Big) < g^{n_x}\Big(h\big(f^{-n_x}(x)\big)\Big) < g^{n_x}(a)\\
\Longrightarrow0 < m^{n_x+1}a < g^{n_x}\Big(h\big(f^{-n_x}(x)\big)\Big) < m^{n_x}a= m^{n_x-1}(ma)
\end{eqnarray*}
That is, $h(x)\in\Big(0,g(a)\Big)$. Also, the function $h$ is a homeomorphism because it is a composition of the homeomorphisms.\\
So the function $h$ was defined on $(0,+\infty)$. We define $h(x)=-h(-x)$ for each $x\in(-\infty,0]$.\\
Now we show that $hof=goh$.\\Suppose $x>a$ is arbitrary. We showed that there exists $n_x\in\Bbb N$ such that
$f^{n_x}(x)\in [f(a),a]$; that is, $f(a)\leq k^{n_x}x\leq a$, so $f(a)\leq k^{n_x-1}(kx)\leq a$; that is, $f(a)\leq k^{n_x-1}(f(x))\leq a$, meaning that $f^{n_x-1}\Big(f(x)\Big) \in [f(a),a]$. We claim that $n_{f(x)}=n_x-1$. We prove this claim with the demonstration by contradiction. Assume there exists natural number $m<n_x-1$ such that $f^{m}\Big(f(x)\Big)\in [f(a), a]$; that is, $f(a)\leq f^{m+1}(x)\leq a$; and know $m+1<n_x$, but this is contradictory with the smallest number of $n_x$ for $x$ and so the claim was proved. Therefore for $x>a$ we have:
\begin{eqnarray*}
\begin{array}{ll}
h\Big(f(x)\Big) & = g^{-n_{f(x)}}\bigg(h\Big(f^{n_{f(x)}}\big(f(x)\big)\Big)\bigg) = g^{-n_x+1}\bigg(h\Big(f^{n_x-1}\big(f(x)\big)\Big)\bigg)\\
   & = g\bigg(g^{-n_x}\Big(h\big(f^{n_x}(x)\big)\Big)\bigg)= g\Big(h(x)\Big).
\end{array}
\end{eqnarray*}
Now let $x\in\Big(o, f(a)\Big)$. We showed that there exists $n_x\in\Bbb N$ such that $f^{-n_x}(x) \in [f(a), a]$; that is, $f(a)\leq f^{-n_x}(x)\leq a$, so $f(a)\leq f^{-(n_x+1)}\big(f(x)\big)\leq a$, We claim that $n_{f(x)}=n_x+1$. We prove this claim with the demonstration by contradiction. Assume that there exists natural number $m < n_x+1$ such that
$f(a)\leq f^{-m}\Big(f(x)\Big)\leq a$ then $f(a)\leq f^{-(m-1)}(x)\leq a$ and have $m-1 < n_x$, but this is contradictory with the way of choice $n_x$ for $x$ and so the claim is proved. Therefore
for $x\in\Big(o, f(a)\Big)$ we have
\begin{eqnarray*}
\begin{array}{ll}
h\Big(f(x)\Big) & = g^{n_{f(x)}}\bigg(h\Big(f^{-n_{f(x)}}\big(f(x)\big)\Big)\bigg)= g^{n_x+1}\bigg(h\Big(f^{-n_x-1}\big(f(x)\big)\Big)\bigg)\\
    & = g^{n_x+1}\bigg(h\Big(f^{-n_x}(x)\Big)\bigg)= g\bigg(g^{n_x}\Big(h\big(f^{-n_x}(x)\big)\Big)\bigg)= g\Big(h(x)\Big).
\end{array}
\end{eqnarray*}
By considering the criteria of functions $f$, $g$ and $h$, we observe that these functions are odd functions. Thus for each $x\in(-\infty, 0]$ we obtain
$h\Big(f(x)\Big)= h\Big(-f(-x)\Big)= -h\Big(f(-x)\Big)= -g\Big(h(-x)\Big)= g\Big(-h(-x)\Big)= g\Big(h(x)\Big).$
Hence, we found the homeomorphism $h$ from $\Bbb R$ to $\Bbb R$ such that $hof=goh$; that is, $f$ and $g$ are topologically conjugate. $\Box$
\end{pf}
By considering the previous lemma, we can say that the expansive functions on $\Bbb{R}$ are topologically conjugate.
\begin{cor}\label{corfunctions}
Suppose
$\left\{ \begin{array}{rl}
 f:\Bbb R & {\longrightarrow \Bbb R}\\
 x & \longmapsto kx
\end{array}\right.$ and
$\left\{ \begin{array}{rl}
 g:\Bbb R & {\longrightarrow \Bbb R}\\
 x & \longmapsto mx
\end{array}\right.$, where $k,m>1$. Then  $f$ and $g$ are topologically conjugate.
\end{cor}
\begin{pf}
Clearly, functions $f$ and $g$ are invertible; we have
$\left\{ \begin{array}{rl}
 f^{-1}:\Bbb R & {\longrightarrow \Bbb R}\\
 x & \longmapsto {1\over k}x
\end{array}\right.$ and
$\left\{ \begin{array}{rl}
 g^{-1}:\Bbb R & {\longrightarrow \Bbb R}\\
 x & \longmapsto {1\over m}x
\end{array}\right.$ such that $0 < {1\over k}, {1\over m} <1$. So by using Lemma\ref{functions}, the functions $f^{-1}$ and $g^{-1}$ are topologically conjugate; that is, there exists a homeomorphism $h$ from $\Bbb R$ to $\Bbb R$ such that $hof^{-1}=g^{-1}oh$, this implies that $f^{-1}= h^{-1}og^{-1}oh$. Therefore, we obtain $foh^{-1}=h^{-1}og$ and since $h^{-1}$ is a homeomorphism hence, $f$ and $g$ are topologically conjugate. $\Box$
\end{pf}
We saw that if $k$ and $m$ both of them are belong to the interval $(0, 1)$ or $(1, +\infty)$ then $f$ and $g$ are topologically conjugate. Moreover similarly, this statement are proved when $k$ and $m$ are belong to the interval $(-1, 0)$ or $(-\infty, -1)$.
Notice that if $k$ and $m$ do not belong to the same interval then $f$ and $g$ are not topologically conjugate. We introduce some examples to show this matter.
\begin{example}
Consider $f(x)= 2x$ and $g(x)= {1\over 2 }x$. Suppose $f$ and $g$ are topologically conjugate, thus those have the same behavior but notice that we have
$\lim_{n\rightarrow\infty}f^{n}(x)= \lim_{n\rightarrow\infty}2^{n}x= \pm\infty$ and
$\lim_{n\rightarrow\infty}g^{n}(x)= \lim_{n\rightarrow\infty}({1\over 2})^{n}x= 0.$
Hence, we may conclude that $f$ and $g$ are not topologically conjugate.
\end{example}
\begin{example}
Assume $f(x)= 3x$ and $g(x)= -3x$. As for the criterion of the function $f$, this function keeps direction but the function $g$ reverse the direction, it means that $f$ and $g$ do not have the same behavior thus those are not topologically conjugate.
\end{example}
\begin{example}
Consider $f(x)= -4x$ and $g(x)= -{1\over 4 }x$. Clearly, we have
$\lim_{n\rightarrow\infty}f^{n}(x)= \lim_{n\rightarrow\infty}({-4})^{n}x= \pm\infty$ and
$\lim_{n\rightarrow\infty}g^{n}(x)= \lim_{n\rightarrow\infty}(-{1\over 4})^{n}x= 0.$
That is, $f$ and $g$ do not have the same behavior thus those can not be topologically conjugate.
\end{example}
\begin{example}
Assume $f(x)= {1\over 5}x$ and $g(x)= -{1\over 5}x$. As regards the function $f$ keeps direction but the function $g$ reverse the direction, so those do not have the same behavior thus we deduce that $f$ and $g$ are not topologically conjugate.
\end{example}

\section{Essential Definitions and Theorems about topological conjugacy of IFSs}
Now we define the concepts of IFS and contractive IFS accurately and formally.(see\cite{Nikiel2007iterated})
\begin{dfn}
Let $(M,d)$ be a complete metric space and $\mathcal{F}$ be a family of continuous mapping $f_{\lambda}:M\rightarrow M$ for every $\lambda \in \Lambda$, where $\Lambda$ is a finite nonempty set; that is, $\mathcal{F}=\Big\{f_{\lambda},\,M\,: \lambda\in\Lambda=\{1,2,\ldots,N\}\Big\}$. We call this family an $\emph{\textbf{Iterated Function System}}$ or shortly, IFS.
\end{dfn}
\begin{dfn}
IFS $\mathcal{F}=\Big\{f_{\lambda},\,M\,: \lambda\in\Lambda\Big\}$ is called $\emph{\textbf{contractive}}$ if each the function $f_{\lambda}$, $\lambda\in\Lambda$, be a contractive function. That is, there exists a positive real number $0<s_{\lambda}<1$ such that for every $x,y\in M$, $d\Big(f_{\lambda}(x),f_{\lambda}(y)\Big)\leq s_{\lambda} d(x,y).$
\end{dfn}
Let $T=\Bbb Z$ or $T=\Bbb N$. $\Lambda^{T}$ denote the set of all infinite sequences $\{\lambda_i\}_{i\in T}$ that $\lambda_i$ is an arbitrary element of $\Lambda$. If $T=\Bbb N$ then every element $\Lambda^{\Bbb N}$ can be showed as $\sigma=\{\lambda_1, \lambda_2,\ldots\}$. Also our intent of the notation $F_{\sigma_n}$ is $F_{\sigma_n}= f_{\lambda_n}o f_{\lambda_{n-1}}o\ldots o f_{\lambda_2}o f_{\lambda_1}$ for every $n\in\Bbb N$.\\
In this paper, we are going to define the concept of topological conjugacy for the IFSs. Previously, this concept has been defined in \cite{Fatehi2015}, but we give comprehensive definition that it includes the previous definition and so we name it as weakly  topological conjugate. The previous definition is as follows:
\begin{dfn}
Suppose $\mathcal{F}=\Big\{f_{\lambda},\,M\,: \lambda\in\Lambda\Big\}$ and $\mathcal{G}=\Big\{g_{\lambda},\,M\,: \lambda\in\Lambda\Big\}$ be two IFSs. The IFSs $\mathcal{F}$ and $\mathcal{G}$ be said topologically conjugate if there exists a homeomorphism $h\,:M\rightarrow M$ such that $f_{\lambda}oh=hog_{\lambda}$ for every $\lambda\in\Lambda$.
\end{dfn}
Our definition is as follows:
\begin{dfn}
Suppose $\mathcal{F}=\Big\{f_{\lambda},\,M\,: \lambda\in\Lambda\Big\}$ and $\mathcal{G}=\Big\{g_{\lambda},\,M\,: \lambda\in\Lambda\Big\}$ are two IFSs. For given $\sigma\in\Lambda^{\Bbb{N}}$, we say that $\mathcal{F}$ and $\mathcal{G}$ are $\emph{\textbf{weakly topological conjugate}}$ if for every $n\in\Bbb N$ there is a homeomorphism $h:M \rightarrow M$ such that $hoF_{\sigma_n}= G_{\sigma_n}oh$.
\end{dfn}
A comparison of two definitions shows that if any two IFSs be topologically conjugate then they will be weakly  topological conjugates. The main problem of the first definition is the presentation of a homeomorphism $h$ for all $\lambda\in\Lambda$, that is a very hard task. We solve this problem by providing a new definition.\\
Hereafter, we will investigate IFSs. In the following, we show that if the model of every two IFSs are $\{ax, bx,\, \Bbb R\}$ that $a$ and $b$ both of them for two IFSs belong to the same interval $(0, 1)\Big(or (-1, 0) or (1, +\infty) or (-\infty, -1)\Big)$ then they are weakly topological conjugates.
\begin{thm}\label{ConIFS}
Suppose $\mathcal{F}= \{k_1x, k_2x,\,\Bbb R\}$ and $\mathcal{G}= \{m_1x, m_2x,\,\Bbb R\}$ are two IFSs where $0<k_i, m_i<1$, $i=1, 2$. Then $\mathcal{F}$ and $\mathcal{G}$ are weakly topological conjugates.
\end{thm}
\begin{pf}
Put $f_i(x)=k_ix$ and $g_i(x)=m_ix$ for $i=1, 2$. Assume $\sigma=\{\lambda_1, \lambda_2,\ldots\}$ is an arbitrary  sequence from indices $\Lambda= \{1, 2\}$. Let $n\in\Bbb{N}$. We know $F_{\sigma_n} = f_{\lambda_n}of_{\lambda_{n-1}}o\ldots of_{\lambda_2}of_{\lambda_1}$, so for every $x\in\Bbb R$ we have
\begin{eqnarray*}
\begin{array}{ll}
F_{\sigma_n}(x) & = f_{\lambda_n}of_{\lambda_{n-1}}o\ldots of_{\lambda_2}of_{\lambda_1}(x)= f_{\lambda_n}of_{\lambda_{n-1}}o\ldots of_{\lambda_2}\Big(f_{\lambda_1}(x)\Big)\\
                & = f_{\lambda_n}of_{\lambda_{n-1}}o\ldots of_{\lambda_2}\Big(k_{\lambda_1}x\Big)= f_{\lambda_n}of_{\lambda_{n-1}}o\ldots of_{\lambda_3}\Big(k_{\lambda_2}.k_{\lambda_1}x\Big)\\
                & = \ldots= k_{\lambda_n}.k_{\lambda_{n-1}}.\ldots .k_{\lambda_2}.k_{\lambda_1}x.
\end{array}
\end{eqnarray*}
Put $k_{\sigma_n}^{*}= k_{\lambda_n}.k_{\lambda_{n-1}}.\ldots .k_{\lambda_2}.k_{\lambda_1}$. Clearly $0<k_{\sigma_n}^{*}<1$ since every $k_{\lambda_{i}}$, $i=1, 2,\ldots,n$, is the value between zero and one.
Also, on the same way for the IFS $\mathcal{G}$ we obtain
\begin{eqnarray*}
\begin{array}{ll}
G_{\sigma_n}(x) & = g_{\lambda_n}o g_{\lambda_{n-1}}o\ldots o g_{\lambda_2}o g_{\lambda_1}(x)= g_{\lambda_n}o g_{\lambda_{n-1}}o\ldots o g_{\lambda_2}\Big(g_{\lambda_1}(x)\Big)\\
                & = g_{\lambda_n}o g_{\lambda_{n-1}}o\ldots o g_{\lambda_2}\Big(m_{\lambda_1}x\Big)= g_{\lambda_n}o g_{\lambda_{n-1}}o\ldots o g_{\lambda_3}\Big(m_{\lambda_2}.m_{\lambda_1}x\Big)\\
                & = \ldots = m_{\lambda_n}.m_{\lambda_{n-1}}.\ldots .m_{\lambda_2}.m_{\lambda_1}x.
\end{array}
\end{eqnarray*}
Now we set $m_{\sigma_n}^{*}= m_{\lambda_n}.m_{\lambda_{n-1}}.\ldots .m_{\lambda_2}.m_{\lambda_1}$. Clearly $0<m_{\sigma_n}^{*}<1$ since every $\lambda_{i}$, $i=1, 2,\ldots,n$, is the value between zero and one.\\
Hence $F_{\sigma_n}(x)=k_{\sigma_n}^{*}x$ and $G_{\sigma_n}(x)=m_{\sigma_n}^{*}x$ where $0<k_{\sigma_n}^{*}, m_{\sigma_n}^{*}<1$, so for every $n\in\Bbb N$ by using Lemma\ref{functions} the functions $F_{\sigma_n}$ and $G_{\sigma_n}$ are topological conjugates; it means that for every $n\in\Bbb N$ there exists a homeomorphism $h$ from $\Bbb R$ to  $\Bbb R$ such that $hoF_{\sigma_n}=G_{\sigma_n}oh$ and this shows that two IFSs, $\mathcal{F}$ and $\mathcal{G}$, are weakly topological conjugates. $\Box$
\end{pf}
\begin{thm}\label{NConIFS}
Suppose $\mathcal{F}= \{-k_1x, -k_2x,\,\Bbb R\}$ and $\mathcal{G}= \{-m_1x, -m_2x,\,\Bbb R\}$ are two IFSs where $0<k_i, m_i<1$, $i=1, 2$. Then $\mathcal{F}$ and $\mathcal{G}$ are weakly topological conjugates.
\end{thm}
\begin{pf}
Put $f_i(x)=-k_ix$ and $g_i(x)=-m_ix$ for $i=1, 2$. Assume $\sigma=\{\lambda_1, \lambda_2,\ldots\}$ is an arbitrary  sequence from indices $\Lambda= \{1, 2\}$. Analogous of the proof of Theorem\ref{ConIFS} for every $n\in\Bbb N$ we obtain
$F_{\sigma_n}(x) = (-1)^{n}k_{\lambda_n}.k_{\lambda_{n-1}}.\ldots .k_{\lambda_2}.k_{\lambda_1}x$ and
$G_{\sigma_n}(x) = (-1)^{n}m_{\lambda_n}.m_{\lambda_{n-1}}.\ldots .m_{\lambda_2}.m_{\lambda_1}x$, for all $x\in\Bbb R$, where $0< k_{\lambda_i},\,m_{\lambda_i}<1$ for $i=1, 2,\ldots,n$.
Put $\mathcal{F}^{*}= \{k_1x, k_2x,\,\Bbb R\}$ and $\mathcal{G}^{*}= \{m_1x, m_2x,\,\Bbb R\}$. In the basis of Theorem\ref{ConIFS}, for given $\sigma$ at above and for every $n\in\Bbb N$ there exists a homeomorphism $h^{*}$ on $\Bbb R$ such that $h^{*}oF_{\sigma_n}^{*}=G_{\sigma_n}^{*}oh^{*}$. Put $h= -h^{*}$. We claim that $hoF_{\sigma_n}=G_{\sigma_n}oh$.\\
First, notice that for each $x\in\Bbb R$, we have $F_{\sigma_n}(x)= (-1)^{n}F_{\sigma_n}^{*}(x)$ and $G_{\sigma_n}(x)= (-1)^{n}G_{\sigma_n}^{*}(x)$ and also the homeomorphism $h^{*}$ is an odd function of lemma\ref{functions}. Moreover, the functions $F_{\sigma_n}^{*}$ and $G_{\sigma_n}^{*}$ are odd, clearly. Suppose $x\in\Bbb R$ be arbitrary. We prove the claim for two states, when $n$ number is odd and when $n$ number is even.
If $n$ is an odd number then we have:
\begin{eqnarray*}
\begin{array}{ll}
h\Big(F_{\sigma_n}(x)\Big) & = h\Big(-F_{\sigma_n}^{*}(x)\Big)= -h^{*}\Big(-F_{\sigma_n}^{*}(x)\Big)= h^{*}\Big(F_{\sigma_n}^{*}(x)\Big)= G_{\sigma_n}^{*}\Big(h^{*}(x)\Big)\\
                           & =  G_{\sigma_n}^{*}\Big(-h^{*}(-x)\Big)= -G_{\sigma_n}^{*}\Big(h^{*}(-x)\Big)=  -G_{\sigma_n}^{*}\Big(-h^{*}(x)\Big)=  G_{\sigma_n}\Big(h(x)\Big).
\end{array}
\end{eqnarray*}
Similarly, when $n$ is an even number, we have:
$$h\Big(F_{\sigma_n}(x)\Big)= h\Big(F_{\sigma_n}^{*}(x)\Big)= -h^{*}\Big(F_{\sigma_n}^{*}(x)\Big)= -G_{\sigma_n}^{*}\Big(h^{*}(x)\Big)=  G_{\sigma_n}^{*}\Big(-h^{*}(x)\Big)=  G_{\sigma_n}\Big(h(x)\Big).$$
So, for every $n\in\Bbb N$ we found the homeomorphism $h$( since $h^{*}$ is the homeomorphism) such that $hoF_{\sigma_n}=G_{\sigma_n}oh$; that is, $\mathcal{F}$ and $\mathcal{G}$ are weakly topological conjugates. $\Box$
\end{pf}
\begin{cor}\label{corExtension}
Suppose $\mathcal{F}= \{k_1x, k_2x,\,\Bbb R\}$ and $\mathcal{G}= \{m_1x, m_2x,\,\Bbb R\}$ where  $k_i$ and $m_i$ are  more than 1 or both of them are less than -1 for each $i=1, 2$. Then $\mathcal{F}$ and $\mathcal{G}$ are weakly topological conjugates.
\end{cor}
\begin{pf}
First, we suppose that $k_i$ and $m_i$  are  more than 1 for each $i=1, 2$. Put $f_i(x)=k_ix$ and $g_i(x)=m_ix$ for $i=1, 2$. Assume $\sigma=\{\lambda_1, \lambda_2,\ldots\}$ is an arbitrary  sequence from indices $\Lambda= \{1, 2\}$. Analogous of the proof of Theorem\ref{ConIFS} for every $n$, we obtain $F_{\sigma_n}(x)=k_{\sigma_n}^{*}x$ where $k_{\sigma_n}^{*}= k_{\lambda_n}.\ldots .k_{\lambda_1}$ that clearly $k_{\sigma_n}^{*}>1$ and also $G_{\sigma_n}(x)=m_{\sigma_n}^{*}x$ where $m_{\sigma_n}^{*}= m_{\lambda_n}.\ldots .m_{\lambda_1}$ that clearly $m_{\sigma_n}^{*}>1$. Thus for every $n\in\Bbb N$, $F_{\sigma_n}$ and $G_{\sigma_n}$ are topologically conjugate from Corollary\ref{corfunctions}; that is, for every $n\in\Bbb N$ there exists homeomorphism $h$ on $\Bbb R$ such that $hoF_{\sigma_n}=G_{\sigma_n}oh$, so $\mathcal{F}$ and $\mathcal{G}$ are weakly topological conjugates.
Now, suppose $k_i$ and $m_i$ are less than -1 for each $i=1, 2$. By considering the previous case we can prove it the similar to the proof of Theorem\ref{NConIFS}. $\Box$
\end{pf}

\section{Extension of Hartman-Grobman Theorem for IFSs}

In the literature has been showed that nonlinear systems sometimes "look like" their linearizations near hyperbolic fixed point( for example in \cite{Palis1980,PERKO,Meiss2008}),this theorem is well known as Hartman-Grobman theorem.\\
$\emph{\textbf{Hartman-Grobman Theorem}}$\,\cite{Guckenheimer1983}\\
 Suppose $x _0$ is hyperbolic fixed point of local $C^{1}$ diffeomorphism $f$ defined on a neighborhood $U$ of $x_o$ in $\Bbb {R}^{m}$. Let $L= Df(x_0)$. Then there exists a neighborhood $U_1\subseteq U$ of $x_0$ and a homeomorphism $h$ from $U_1$ into $\Bbb {R}^{m}$ such that $h(x_0)= 0$ and $hf(x)= Lh(x)$ for $x\in U_1\cap f^{-1}(U_1)$.$\Big($ or $hfh^{-1}(y)= L(y)$ for $h^{-1}(y)\in U_1\cap f^{-1}(U_1)$.$\Big)$\\

Notice down theorems and their corollaries in order to extend Hartman-Grobman Theorem for IFSs.
\begin{thm}\label{mainTh1}
Suppose $\mathcal{F}= \{k_1I+ \varphi_1, k_2I+ \varphi_2; \Bbb R\}$ and $\mathcal{G}= \{m_1I+ \psi_1, m_2I+ \psi; \Bbb R\}$ are two IFSs where $I$ is identity map on $\Bbb R$ and for $i=1, 2$, $k_i$ and $m_i$ all of them have the same sign and $0<\mid k_i\mid, \mid m_i\mid<1$ and also the functions $\varphi_i$ and $\psi_i$, $i=1, 2$, are Lipschitz functions with Lipschitz constant at most $\epsilon$ that those contain the conditions
$\varphi_i(0)= \psi_i(0)= 0$ and $0< \mid k_i\mid+ \epsilon, \mid m_i\mid+ \epsilon< 1$ for each $i=1, 2$.\\Then
\begin{enumerate}
\item  the functions $k_iI+ \varphi_i$ and $m_iI+ \psi_i$ are contractions for $i=1, 2$,\\
\item $\mathcal{F}$ and $\mathcal{G}$ are weakly topological conjugates.
\end{enumerate}
\end{thm}
\begin{pf}
{\bf 1.} Consider the usual norm $\parallel.\parallel$ on $\Bbb R$. Since the functions $\varphi_i$ and $\psi_i$, $i=1, 2$, are Lipschitz thus for every $x,y\in\Bbb R$ we have $\parallel\varphi_i(x)- \varphi_i(y)\parallel < \epsilon\parallel x- y\parallel$ and $\parallel\psi_i(x)- \psi_i(y)\parallel < \epsilon\parallel x- y\parallel$ for $i=1, 2$. Then
\begin{eqnarray*}
\begin{array}{ll}
\parallel\Big(k_i I+ \varphi_i\Big)(x)- \Big(k_i I+ \varphi_i\Big)(y)\parallel & = \parallel k_i x+ \varphi_i(x)- k_i y- \varphi_i(y)\parallel\\
                           & = \parallel k_i (x-y)+ \varphi_i(x)- \varphi_i(y)\parallel\\
                           & \leq\parallel k_i (x-y)\parallel+ \parallel\varphi_i(x)- \varphi_i(y)\parallel \\
                          & <\mid k_i\mid\parallel x-y \parallel+ \epsilon\parallel x-y \parallel
                           =\Big(\mid k_i\mid+ \epsilon\Big)\parallel x-y \parallel
\end{array}
\end{eqnarray*}
Therefore, by considering the hypothesis of the theorem; that is, $0< k_i+ \epsilon < 1$ for $i=1, 2$, the previous relation shows that the function $k_i I+ \varphi_i$ is a contraction and similarly we obtain that the function $m_i I+ \psi_i$ is a contraction for each $i=1, 2$ and so the first statement is proved.\\
{\bf 2.} Assume that $\sigma=\{\lambda_1, \lambda_2,\ldots\}$ is an arbitrary sequence from indices $\Lambda= \{1, 2\}$. First, we show that ${\{\mid F_{\sigma_n} \mid\}}_{n=1}^{\infty}$ is a strictly decreasing sequence and the sequence
${\{F_{\sigma_n} \}}_{n=1}^{\infty}$ is convergence to zero.\\
Suppose $x\in\Bbb R$ is arbitrary. Now we write some of terms of the sequence ${\{\mid F_{\sigma_n} \mid\}}_{n=1}^{\infty}$;
\begin{eqnarray*}
\begin{array}{ll}
\mid F_{\sigma_1}(x)\mid & = \mid f_{\lambda_1}(x)\mid= \mid k_{\lambda_1}x+ \varphi_{\lambda_1}(x)\mid\\
\mid F_{\sigma_2}(x)\mid & = \mid f_{\lambda_2}\Big(f_{\lambda_1}(x)\Big)\mid= \mid k_{\lambda_2}\Big(k_{\lambda_1}x+ \varphi_{\lambda_1}(x)\Big)+ \varphi_{\lambda_2}\Big(k_{\lambda_1}x+ \varphi_{\lambda_1}(x)\Big)\mid\\
                         & \leq\mid k_{\lambda_2}\Big(k_{\lambda_1}x+ \varphi_{\lambda_1}(x)\Big)\mid+ \mid\varphi_{\lambda_2}\Big(k_{\lambda_1}x+ \varphi_{\lambda_1}(x)\Big)\mid.
\end{array}
\end{eqnarray*}
By using the suppositions of the theorem; that is, for $i=1, 2$, $\varphi_i(0)= 0$ and $0< \mid k_i\mid+ \epsilon< 1$ and also the functions $\varphi_i$ are Lipschitz with constant at most $\epsilon$, we can write the previous relation as follows;
\begin{eqnarray*}
\begin{array}{ll}
\mid F_{\sigma_2}(x)\mid & < \mid k_{\lambda_2}\mid\mid k_{\lambda_1}x+ \varphi_{\lambda_1}(x)\mid+ \epsilon \mid k_{\lambda_1}x+ \varphi_{\lambda_1}(x)\mid\\
                         & = \Big(\mid k_{\lambda_2}\mid+ \epsilon\Big)\mid k_{\lambda_1}x+ \varphi_{\lambda_1}(x)\mid < \mid k_{\lambda_1}x+ \varphi_{\lambda_1}(x)\mid= \mid F_{\sigma_1}(x)\mid.
\end{array}
\end{eqnarray*}
Generally, for every $n$ we have:
\begin{eqnarray*}
\begin{array}{ll}
F_{\sigma_n}(x) & = f_{\lambda_n}of_{\lambda_{n-1}}o\ldots of_{\lambda_2}of_{\lambda_1}(x)
                = f_{\lambda_n}\Big(f_{\lambda_{n-1}}o\ldots of_{\lambda_2}of_{\lambda_1}(x)\Big)\\
                &= f_{\lambda_n}\Big(F_{\sigma_{n-1}}(x)\Big)
                 = k_{\lambda_n}\Big(F_{\sigma_{n-1}}(x)\Big)+ \varphi_{\lambda_n}\Big(F_{\sigma_{n-1}}(x)\Big).
\end{array}
\end{eqnarray*}
So
\begin{eqnarray*}
\begin{array}{ll}
\mid F_{\sigma_n}(x)\mid & = \mid k_{\lambda_n}\Big(F_{\sigma_{n-1}}(x)\Big)+ \varphi_{\lambda_n}\Big(F_{\sigma_{n-1}}(x)\Big)\mid\\
                        & \leq \mid k_{\lambda_n}\mid\mid F_{\sigma_{n-1}}(x)\mid+ \mid \varphi_{\lambda_n}\Big(F_{\sigma_{n-1}}(x)\Big)\mid
                         < \mid k_{\lambda_n}\mid\mid F_{\sigma_{n-1}}(x)\mid+ \epsilon\mid F_{\sigma_{n-1}}(x)\mid\\
                         &= \Big(\mid k_{\lambda_n}+ \epsilon\mid\Big)\mid F_{\sigma_{n-1}}(x)\mid
                         < \mid F_{\sigma_{n-1}}(x)\mid.
\end{array}
\end{eqnarray*}
Hence, the sequence ${\{\mid F_{\sigma_n} \mid\}}_{n=1}^{\infty}$ is a strictly decreasing and since it is bounded from below( for every $x\in\Bbb R$, $\mid F_{\sigma_n}(x)\mid >0$) thus it is convergence.\\
Put $k= Max\{\mid k_1\mid+ \epsilon, \mid k_2\mid+ \epsilon\}$, clearly $0<k<1$. Now, for each $i=1, 2$ we will obtain
$$\parallel\Big(k_i I+ \varphi_i\Big)(x)- \Big(k_i I+ \varphi_i\Big)(y)\parallel < k\parallel x-y \parallel.$$
Also $\Big(k_i I+ \varphi_i\Big)(0)= 0$, thus we have:
\begin{eqnarray*}
\begin{array}{ll}
\mid F_{\sigma_n}(x)\mid & = \mid k_{\lambda_n}\Big(F_{\sigma_{n-1}}(x)\Big)+ \varphi_{\lambda_n}\Big(F_{\sigma_{n-1}}(x)\Big)\mid\\
                         & < k\mid F_{\sigma_{n-1}}(x)\mid = k\mid k_{\lambda_{n-1}}\Big(F_{\sigma_{n-2}}(x)\Big)+ \varphi_{\lambda_{n-1}}\Big(F_{\sigma_{n-2}}(x)\Big)\mid\\
                         & < k.k\mid F_{\sigma_{n-2}}(x)\mid = k^{2}\mid F_{\sigma_{n-2}}(x) \mid.
\end{array}
\end{eqnarray*}
Keep on this way, finally we will obtain:
$$\mid F_{\sigma_n}(x)\mid < k^{n-1}\mid F_{\sigma_1}(x)\mid = k^{n-1}\mid k_{\lambda_1}x+ \varphi_{\lambda_1}(x)\mid < k^{n-1}.k\mid x \mid = k^{n}\mid x \mid.$$
Hence, for every $x\in\Bbb R$,\,$\mid F_{\sigma_n}(x)\mid < k^{n}\mid x\mid$ where $0<k^{n}<1$ and since the functions of $\mathcal{F}$ are contraction, by  part(1) of the theorem, so it implies that the sequence ${\{F_{\sigma_n} \}}_{n=1}^{\infty}$ is convergence to zero. Similarly, the sequence ${\{G_{\sigma_n} \}}_{n=1}^{\infty}$ is convergence to zero. Thus, these two IFSs have the same behavior as two qualified IFSs at Theorem\ref{ConIFS}, consequently $\mathcal{F}$ and $\mathcal{G}$ are weakly topological conjugates. $\Box$
\end{pf}
%\begin{dfn}
%Let $f$ be a $C^{1}$ diffeomorphism from a neighborhood $U$ of $x_0$ in $\Bbb {R}^{n}$ into $\Bbb {R}^{n}$. The point
%$x _0$ is called $\emph{\textbf{fixed point}}$ for $f$ if $f(x_0)=x_0$.
%\end{dfn}
\begin{cor}\label{maincor1}
Suppose $\mathcal{F}=\{f_1, f_2,\,\Bbb R\}$ is an IFS where the functions $f_1$ and $f_2$ are diffeomorphisms on $\Bbb R$. The origin is a fixed point of the functions $f_1$ and $f_2$ and also assume the derivative values of these functions at the origin have the same sign and $0< \mid\acute{f_1}(0)\mid, \mid\acute{f_2}(0)\mid <1$. Consider IFS $\mathcal{G}=\{\acute{f_1}(0)I, \acute{f_2}(0)I,\,\Bbb R\}$. Then $\mathcal{F}$ and $\mathcal{G}$ are weakly topological conjugates on a neighborhood of zero.
\end{cor}
\begin{pf}
Suppose $\epsilon>o$ is a number such that $0<\epsilon + \mid\acute{ f_i}(0)\mid<1$ for each $i=1, 2$. For every $i=1,2$ by considering the lemma(4.4) in \cite{Palis1980}, for given $\epsilon>o$ there exists a neighborhood $U_i$ of zero and an extension of $f_i|_{U_i}$ to $\Bbb R$ of the form $\acute{f_i}(0)I+ \varphi_i$ where $\varphi_i$ is a bounded continuous map from $\Bbb R$ to $\Bbb R$ that it has Lipschitz constant at most $\epsilon$. Since zero is a fixed point of $f_i$ and the functions $f_i$ and  $\acute{f_i}(0)I+ \varphi_i$ are equal on $U$, therefore $\varphi_i(0)=f_i(0)= 0$. Now, put $U= U_1\cap U_2$ and $\mathcal{F^{*}}=\{\acute{f_1}(0)I+ \varphi_1, \acute{f_2}(0)I+ \varphi_2,\,\Bbb R\}$. On the basis of the previous theorem we conclude that $\mathcal{F^{*}}$ and $\mathcal{G}$ are weakly topological conjugates and since IFS $\mathcal{F}$ has the same behavior as IFS $\mathcal{F^{*}}$ on $U$( because the functions IFS $\mathcal{F^{*}}$ are extended the functions IFS $\mathcal{F}$ on $U$) so IFSs $\mathcal{F}$ and $\mathcal{G}$ are weakly topological conjugates on $U$ and the statement is proved. $\Box$
\end{pf}
\begin{cor}\label{maincor2}
Suppose we have the assumptions of the previous corollary, only $\mid\acute{f_1}(0)\mid$,  $\mid\acute{f_2}(0)\mid>1$.
Consider IFS $\mathcal{G}=\{\acute{f_1}(0)I, \acute{f_2}(0)I,\,\Bbb R\}$. Then $\mathcal{F}$ and $\mathcal{G}$ are weakly topological conjugates on a neighborhood of zero.
\end{cor}
\begin{pf}
Since the functions of IFS $\mathcal{F}$ are diffeomorphisms, so we can consider IFS\,$\mathcal{F}^{*}=\{{f_1}^{-1}, {f_2}^{-1},\,\Bbb R\}$. Clearly, the origin is fixed point of the functions ${f_1}^{-1}$ and ${f_2}^{-1}$. We know $({{f_i}^{-1})}^{'}(0)= {{1} /{\acute{f_i}(0)}}$, therefore the values $({{f_1}^{-1})}^{'}(0)$ and $({{f_2}^{-1})}^{'}(0)$ have the same sign and $0<\mid({{f_1}^{-1})}^{'}(0)\mid, \mid({{f_2}^{-1})}^{'}(0)\mid<1$. Thus $\mathcal{F}^{*}$ and IFS\, $\mathcal{G}^{*}=\{({{f_1}^{-1})}^{'}(0)I, ({{f_2}^{-1})}^{'}(0)I,\,\Bbb R\}$ are weakly topological conjugates on neighborhood $U$ of zero, of the previous corollary; that is, for every $\sigma$ and $n\in\Bbb N$ there exists a homeomorphism $h$ on $U$ such that $hoF_{\sigma_n}^{*}=G_{\sigma_n}^{*}oh$. Now for $n\in\Bbb N$ put $\sigma^{*}= \{\lambda_n, \lambda_{n-1},\ldots,\lambda_1,\lambda_{n+1},\ldots\}$. For these $n$ and $\sigma^{*}$ there exists homeomorphism $h$ on $U$ such that $hoF_{\sigma_n}^{*}=G_{\sigma_n}^{*}oh$. So we have
$F_{\sigma_n}^{*}(x)=  f_{\lambda_1}^{-1}of_{\lambda_2}^{-1}o\ldots of_{\lambda_n}^{-1}(x)= \Big( f_{\lambda_n}of_{\lambda_{n-1}}o\ldots of_{\lambda_1}\Big)^{-1}(x)= F_{\sigma_n}^{-1}(x)$
and $G_{\sigma_n}^{*}(x)= k^{*}x$ where $0<\mid k^{*}\mid<1$ and clearly, $G_{\sigma_n}^{*}(x)= G_{\sigma_n}^{-1}(x)$. Thus we can obtain that ${\Big(G_{\sigma_n}^{*}\Big)}^{-1}oh= ho{\Big(F_{\sigma_n}^{*}\Big)}^{-1}$ and subsequently $hoF_{\sigma_n}=G_{\sigma_n}oh$; that is, $\mathcal{F}$ and $\mathcal{G}$ are weakly topological conjugates on neighborhood $U$ of zero. $\Box$
\end{pf}
\begin{thm}\label{mainTh2}
Assume $\mathcal{F}=\{f_1, f_2,\,\Bbb R\}$ is an IFS where the functions $f_1$ and $f_2$ are homeomorphisms on $\Bbb R$. The origin is a fixed point of the functions $f_1$ and $f_2$. Also suppose $\acute{f_1}(0)$ and $\acute{f_2}(0)$ have the same sign and $0<\mid\acute{f_1}(0)\mid<1$ and $\mid\acute{f_2}(0)\mid>1$. Consider the IFS $\mathcal{G}=\{\acute{f_1}(0)I, \acute{f_2}(0)I,\,\Bbb R\}$. Let $\sigma=\{\lambda_1, \lambda_2,\ldots\}$ and numbers of times of $\lambda_i$,$i\in\Bbb N$, that $\lambda_i= 1$  is $n_1$ and numbers of times of $\lambda_i$ that $\lambda_i= 2$ is $n_2$ such that $‎\lim_{n\rightarrow +\infty}{{n_1}/{n_2}}= +\infty$ (or $‎\lim_{n\rightarrow +\infty}{{n_2}/{n_1}}= o$). Then for every $n\in\Bbb N$ there exists a homeomorphism $h$ on a neighborhood of zero such that $hoF_{\sigma_n}= G_{\sigma_n}oh$.
\end{thm}
\begin{pf}
Put $\acute{f_i}(0)= a_i$,$i=1, 2$. Assume $\epsilon>o$ is a number such that $\mid a_1\mid+\epsilon <1$
and $\mid a_2\mid-\epsilon >1$. Let the neighborhood $U$, the functions $\varphi_1$ and $\varphi_2$ and also IFS $\mathcal{F^{*}}$ be whose were introduced at Corollary\ref{maincor1}. Put $f_i^{*}= a_i I+ \varphi_i$, $i=1, 2$. We prove that for every $x\in\Bbb R$ sequences ${\{F_{\sigma_n}^{*}(x) \}}_{n=1}^{\infty}$ and ${\{G_{\sigma_n}(x) \}}_{n=1}^{\infty}$ are convergent to zero.\\
Suppose $x\in\Bbb R$ is arbitrary, we write some of terms the sequence ${\{F_{\sigma_n}^{*}(x) \}}_{n=1}^{\infty}$:
\begin{eqnarray*}
\begin{array}{ll}
\mid F_{\sigma_1}^{*}(x)\mid & = \mid f_{\lambda_1}^{*}(x)\mid= \mid a_{\lambda_1}x+ \varphi_{\lambda_1}(x)\mid\leq \mid a_{\lambda_1}\mid\mid x\mid+ \epsilon\mid x\mid= \Big(\mid a_{\lambda_1}\mid+ \epsilon\Big)\mid x\mid\\
\mid F_{\sigma_2}^{*}(x)\mid & = \mid f_{\lambda_2}^{*}\Big(f_{\lambda_1}^{*}(x)\Big)\mid =
\mid a_{\lambda_2}\Big(f_{\lambda_1}^{*}(x)\Big)+ \varphi_{\lambda_2}\Big(f_{\lambda_1}^{*}(x)\Big)\mid\\
                       & \leq \mid a_{\lambda_2}\mid\mid f_{\lambda_1}^{*}(x)\mid+ \epsilon\mid f_{\lambda_1}^{*}(x)\mid=                        \Big(\mid a_{\lambda_2}\mid+ \epsilon\Big)\mid f_{\lambda_1}^{*}(x)\mid\\
                        & \leq \Big(\mid a_{\lambda_2}\mid+ \epsilon\Big).\Big(\mid a_{\lambda_1}\mid+ \epsilon\Big)\mid x\mid.
\end{array}
\end{eqnarray*}
Applying induction we get
$$\mid F_{\sigma_{n-1}}^{*}(x)\mid\leq\Big(\mid a_{\lambda_{n-1}}\mid+ \epsilon\Big).\Big(\mid a_{\lambda_{n-2}}\mid+ \epsilon\Big).\ldots.\Big(\mid a_{\lambda_1}\mid+ \epsilon\Big)\mid x\mid.$$
So
\begin{eqnarray*}
\begin{array}{ll}
\mid F_{\sigma_n}^{*}(x)\mid & = \mid f_{\lambda_n}^{*}of_{\lambda_{n-1}}^{*}o\ldots of_{\lambda_2}^{*}of_{\lambda_1}^{*}(x)\mid
                 = \mid f_{\lambda_n}^{*}\Big(f_{\lambda_{n-1}}^{*}o\ldots of_{\lambda_2}^{*}of_{\lambda_1}^{*}(x)\Big)\mid\\
                & = \mid f_{\lambda_n}^{*}\Big(F_{\sigma_{n-1}}^{*}(x)\Big)\mid
                 = \mid a_{\lambda_n}\Big(F_{\sigma_{n-1}}^{*}(x)\Big)+ \varphi_{\lambda_n}\Big(F_{\sigma_{n-1}}^{*}(x)\Big)\mid\\
                &\leq\mid a_{\lambda_n}\mid\mid F_{\sigma_{n-1}}^{*}(x)\mid+ \epsilon\mid F_{\sigma_{n-1}}^{*}(x)\mid
                 = \Big(\mid a_{\lambda_n}\mid+ \epsilon\Big)\mid F_{\sigma_{n-1}}^{*}(x)\mid\\
                & \leq\Big(\mid a_{\lambda_n}\mid+ \epsilon\Big).\Big(\mid a_{\lambda_{n-1}}\mid+ \epsilon\Big).\ldots.\Big(\mid a_{\lambda_1}\mid+ \epsilon\Big)\mid x\mid.
\end{array}
\end{eqnarray*}
By utilizing supposition we can write the previous relation as follows:
$$\mid F_{\sigma_n}^{*}(x)\mid\leq{\Big(\mid a_1\mid+ \epsilon\Big)}^{n_1}.{\Big(\mid a_2\mid+ \epsilon\Big)}^{n_2}\mid x\mid$$
In the basis assumption $\lim_{n\rightarrow +\infty}{{n_1}\over {n_2}}= +\infty$, $n_1$ is very larger than $n_2$ when $n\rightarrow +\infty$; that is, $n_1$ gradually approaches $n$, then from the relations $\mid a_1\mid+ \epsilon <1$, $\mid a_2\mid>1$ and
$\mid F_{\sigma_n}^{*}(x)\mid\leq{\Big(\mid a_1\mid+ \epsilon\Big)}^{n_1}.{\Big(\mid a_2\mid+ \epsilon\Big)}^{n-n_1}\mid x\mid$, we conclude that $F_{\sigma_n}^{*}(x)\rightarrow o$ as $n\rightarrow +\infty$.
Also for IFS $\mathcal{G}$ and every $x\in\Bbb R$, we have
\begin{eqnarray*}
\begin{array}{ll}
\mid G_{\sigma_n}(x)\mid & = \mid a_{\lambda_n}.a_{\lambda_{n-1}}.\ldots.a_{\lambda_1}x\mid= \mid a_{\lambda_n}\mid.\mid a_{\lambda_{n-1}}\mid.\ldots.\mid a_{\lambda_1}\mid\mid x\mid\\
                          & = {\mid a_1\mid}^{n_1}{\mid a_2\mid}^{n_2}\mid x\mid= {\mid a_1\mid}^{n_1}{\mid a_2\mid}^{n- n_1}\mid x\mid.
\end{array}
\end{eqnarray*}
We know $0<\mid a_1\mid<1$, so with a reasoning similar to the above argument we obtain $G_{\sigma_n}(x)\rightarrow o$ as $n\rightarrow +\infty$.
Then, these two IFSs have the same behavior as two qualified IFSs at Theorem\ref{ConIFS}, it follow that for given $\sigma\in\Lambda^{\Bbb N}$ and every $n\in\Bbb N$, there exists homeomorphism $h$ such that $hoF_{\sigma_n}^{*}= G_{\sigma_n}oh$. Since IFSs $\mathcal{F^{*}}$ and $\mathcal{F}$ are equal on $U$, thus for given $\sigma$ and every $n\in\Bbb N$ there exists a homeomorphism $h$ on $U$ such that $hoF_{\sigma_n}= G_{\sigma_n}oh$ and therefore the statement is proved. $\Box$
\end{pf}
Notice that hereafter if for a given $\sigma\in\Lambda^{\Bbb N}$ and every $n\in\Bbb N$ there exists a homeomorphism $h$ such that $hoF_{\sigma_n}= G_{\sigma_n}oh$ then we say that IFSs $\mathcal{F}$ and $\mathcal{G}$ are weakly topological conjugates relative to $\sigma$.
\begin{cor}\label{maincor3}
Suppose we have all the assumptions of Theorem\ref{mainTh2} but instead, we have $\lim_{n\rightarrow +\infty}{{n_1}/{n_2}}= 0$. Then $\mathcal{F}$ and $\mathcal{G}$ are weakly topological conjugates relative to $\sigma$.
\end{cor}
\begin{pf}
Using the obtained relations at Theorem\ref{mainTh2}, since $\lim_{n\rightarrow +\infty}{{n_1}/{n_2}}= 0$ we get $F_{\sigma_n}^{*}(x)\rightarrow\infty$ and $G_{\sigma_n}(x)\rightarrow \infty$ then $\mathcal{F}^{*}$ and $\mathcal{G}$ are weakly topological conjugate relative to given $\sigma$ of Corollary\ref{corExtension}. Thus IFSs $\mathcal{F}$ and $\mathcal{G}$ are weakly topological conjugates relative to $\sigma$ on $U$. $\Box$
\end{pf}
\begin{dfn}
Let $f$ be a $C^{1}$ diffeomorphism from a neighborhood $U$ of $x_o$ in $\Bbb {R}^{n}$ into $\Bbb {R}^{n}$. The fixed point $x _0$ is called $\emph{\textbf{hyperbolic}}$ if all the eigenvalues of $Df(x_0)$ have absolute values with norm different from one, \cite{Pilyugin1992introduction}
\end{dfn}
Now, we express generalized Hartman-Grobman Theorem for IFSs.
\begin{thm}\label{main}
$\emph{\textbf{(Generalized Hartman-Grobman Theorem for IFSs)}}$\\
Suppose $\mathcal{F}=\{f_{\lambda}: \lambda\in\Lambda, \Bbb R\}$\,($\Lambda$ is a finite nonempty set) is an IFS and the origin is hyperbolic fixed point of the homeomorphisms $f_\lambda$ for every $\lambda\in\Lambda$. Consider IFS\, $\mathcal{G}=\{\acute{f_\lambda}(0)I\,:\,\lambda\in\Lambda, \Bbb R\}$; we say it "linear part of IFS $\mathcal{F}$". Then
\begin{enumerate}
\item if $\acute{f_\lambda}(0)$ belong to the same interval $(0, 1)\Big($ or $(-1, 0)$ or $(1, +\infty)$ or $(-\infty, -1)\Big)$ for all $\lambda\in\Lambda$, then $\mathcal{F}$ and $\mathcal{G}$ are weakly topological conjugates on a neighborhood of zero.\\
\item suppose $\acute{f_\lambda}(0)$,\,$\lambda\in\Lambda$, all of them have the same sign. Moreover some of them are belong to the same interval $(0, 1)\Big($ or $(-1,0)\Big)$ and some of them are belong to the same interval $(1, +\infty)\Big($ or $(-\infty, -1)\Big)$. Assume that $\sigma\in\Lambda^{\Bbb N}$ be given and numbers of times of $\lambda_i$,$i\in\Bbb N$, that $0< \mid\acute{f_{\lambda_i}}(0)\mid<1$ is $n_1$ and numbers of times of $\lambda_i$ that $\mid\acute{f_{\lambda_i}}(0)\mid >1$ is $n_2$ such that $\lim_{n\rightarrow +\infty}{{n_1}/{n_2}}= +\infty \Big($ or $\lim_{n\rightarrow +\infty}{{n_1}/{n_2}}= o\Big)$ thus $\mathcal{F}$ and $\mathcal{G}$ are weakly topological conjugates relative to $\sigma$ on a neighborhood of zero.
\end{enumerate}
\end{thm}
\begin{pf}
{\bf 1.}\, Similarly, we can extend Theorem\ref{mainTh1} and Corollaries\ref{maincor1} and \ref{maincor2} for the case that $\Lambda$ is a finite nonempty set and subsequently the first statement is true.\\
{\bf 2.}\, With the same way as Theorem\ref{mainTh2} and Corollary\ref{maincor2}, we can see that those also for the case that $\Lambda$ is a finite nonempty set are true, so the second statement is true. $\Box$
\end{pf}
In the following theorem, we examine topological conjugacy for two IFSs.
\begin{thm}\label{mainTh3}
Suppose $\mathcal{F}=\{f_{\lambda}, \Bbb R : \lambda\in\Lambda\}$ and $\mathcal{G}=\{g_{\lambda}, \Bbb R : \lambda\in\Lambda\}$ are two IFSs where for every $\lambda\in\Lambda$ the functions $f_{\lambda}$ and $g_{\lambda}$ are homeomorphisms. Let origin be a fixed point of the functions IFSs  $\mathcal{F}$ and  $\mathcal{G}$. Assume for all $\lambda\in\Lambda$,\, $\acute{f_\lambda}(0)$ and $\acute{g_\lambda}(0)$ are belong to the same interval $(0, 1)\Big($ or $(-1, 0)$ or $(1, +\infty)$ or $(-\infty, -1)\Big)$ then $\mathcal{F}$ and $\mathcal{G}$ are weakly topological conjugates on a neighborhood of zero.
\end{thm}
\begin{pf}
Let $\mathcal{F}^{*}$ be the linear part of the IFS $\mathcal{F}$ and $\mathcal{G}^{*}$ be the linear part of the IFS $\mathcal{G}$. By assumptions and the first part of Generalized Hartman-Grobman Theorem for IFSs, we obtain that IFSs $\mathcal{F}$ and $\mathcal{F}^{*}$ are weakly topological conjugates on neighborhood $U$ of zero and also IFSs $\mathcal{G}$ and $\mathcal{G}^{*}$ are weakly topological conjugates  on neighborhood $V$ of zero. Applying primary theorems of this paper, we can conclude that IFSs $\mathcal{F}^{*}$ and $\mathcal{G}^{*}$ are weakly topological conjugates. Now, put $W=U\cap V$. Clearly $W$ is a neighborhood of zero and thus we get that IFSs $\mathcal{F}$ and $\mathcal{G}$ are weakly topological conjugates on neighborhood $W$ of zero. $\Box$
\end{pf}
\section{Topological Conjugacy of $m$-dimensional IFSs}

Now, we assume that the functions $f_i$ of IFS $\mathcal{F}$ are determined from $\Bbb {R}^{m}$ to $\Bbb {R}^{m}$ and we investigate concept "weakly topological conjugate" for some of special IFSs.
\begin{thm}
Consider IFSs $\mathcal{F}= \{A, B, \Bbb {R}^{m}\}$ and $\mathcal{G}= \{C, D, \Bbb {R}^{m}\}$ where $A, B, C$ and $D$ are diagonal matrix respectively with the diagonal elements $a_{ii}, b_{ii}, c_{ii}$ and $d_{ii}$, $i=1,2,\ldots,m$. If all these elements are belong to the same interval $(0, 1)\Big($ or $(-1, 0)$ or $(1,+\infty)$ or $(-\infty, -1)\Big)$ then $\mathcal{F}$ and $\mathcal{G}$ are weakly topological conjugates.
\end{thm}
\begin{pf}
Assume
$X= \left [\begin{array}{c}
x_1\\
x_2\\
\vdots\\
x_m\\
\end {array}\right ]$
and also $\sigma\in\Lambda^{\Bbb N}$ be arbitrary. Since product of diagonal matrices is a diagonal matrix, so we obtain:\\
 $F_{\sigma_n}(X)= \left [\begin{array}{c}
{a_{11}}^{n_1}{b_{11}}^{n_2}x_1\\
{a_{22}}^{n_1}{b_{22}}^{n_2}x_2\\
\vdots\\
{a_{mm}}^{n_1}{b_{mm}}^{n_2}x_m\\
\end {array}\right ]$
and
$G_{\sigma_n}(X)= \left [\begin{array}{c}
{c_{11}}^{n_1}{d_{11}}^{n_2}x_1\\
{c_{22}}^{n_1}{d_{22}}^{n_2}x_2\\
\vdots\\
{c_{mm}}^{n_1}{d_{mm}}^{n_2}x_m\\
\end {array}\right ]$
where $n_1+ n_2=n$, in fact $n_1$ is numbers of times of iteration $A$(associated to $G_{\sigma_n}(X)$ is $C$) at sequence $\sigma$ and $n_2$ is numbers of times of iteration $B$(associated to $G_{\sigma_n}(X)$ is $D$) at sequence $\sigma$. Now, for $i=1, 2,\ldots,m$ put ${\mathcal{F}}_i=\{a_{ii}x_i I, b_{ii}x_iI, \Bbb R\}$ and ${\mathcal{G}}_i=\{c_{ii}x_iI, d_{ii}x_iI, \Bbb R\}$. On the basis of the previous proved statements, for each $i=1, 2,\ldots,m$ IFSs ${\mathcal{F}}_i$ and ${\mathcal{G}}_i$ are weakly topological conjugates, thus there exists a homeomorphism $h_i:\Bbb R \rightarrow \Bbb R$ such that $h_i\Big(F_{i,\sigma_n}(x_i)\Big)= G_{i,\sigma_n}\Big(h_i(x_i)\Big)$. We define the function $h:\Bbb {R}^{m} \rightarrow \Bbb {R}^{m}$ with criterion
$h(X)= \left [\begin{array}{c}
h_1(x_1)\\
h_2(x_2)\\
\vdots\\
h_m(x_m)\\
\end {array}\right ]$.
Clearly, $h$ is a homeomorphism since for each $i=1, 2,\ldots,m$ the function $h_i$ is the homeomorphism. We claim that the homeomorphism $h$ hold $h\Big(F_{\sigma_n}(X)\Big)= G_{\sigma_n}\Big(h(X)\Big)$ for every $X\in \Bbb {R}^{m}$.\\
For each $X\in \Bbb {R}^{m}$ we have
\begin{eqnarray*}
\begin{array}{ll}
h\Big(F_{\sigma_n}(X)\Big) & = \left [\begin{array}{c}
h_1\Big({a_{11}}^{n_1}{b_{11}}^{n_2}x_1\Big)\\
h_2\Big({a_{22}}^{n_1}{b_{22}}^{n_2}x_2\Big)\\
\vdots\\
h_m\Big({a_{mm}}^{n_1}{b_{mm}}^{n_2}x_m\Big)\\
\end {array}\right ]
= \left [\begin{array}{c}
h_1\Big(F_{1,\sigma_n}(x_1)\Big)\\
h_2\Big(F_{2,\sigma_n}(x_2)\Big)\\
\vdots\\
h_m\Big(F_{m,\sigma_n}(x_m)\Big)\\
\end {array}\right ]\\
                         & = \left [\begin{array}{c}
G_{1,\sigma_n}\Big(h_1(x_1)\Big)\\
G_{2,\sigma_n}\Big(h_2(x_2)\Big)\\
\vdots\\
G_{m,\sigma_n}\Big(h_m(x_m)\Big)
\end {array}\right ]
=  \left [\begin{array}{c}
{c_{11}}^{n_1}{d_{11}}^{n_2}h_1(x_1)\\
{c_{22}}^{n_1}{d_{22}}^{n_2}h_2(x_2)\\
\vdots\\
{c_{mm}}^{n_1}{d_{mm}}^{n_2}h_m(x_m)\\
\end {array}\right ]
= G_{\sigma_n}\Big(h(X)\Big)
\end{array}
\end{eqnarray*}
Therefore the our claim is proved, hence IFSs $\mathcal{F}$ and $\mathcal{G}$ will be weakly topological conjugates. $\Box$
\end{pf}
\begin{thm}
Let $J\subseteq\Bbb N$ be a finite set. Consider the IFS $\mathcal{F}=\{D_j:\,j\in J, \Bbb {R}^{m}\}$ where $D_j$ is a diagonal matrix for every $j\in J$. Let $\mathcal{G}=\{AD_jA^{-1}:\,j\in J, \Bbb {R}^{m}\}$ where matrix $A$ is an invertible matrix. Then $\mathcal{F}$ and $\mathcal{G}$ are weakly topological conjugates.
\end{thm}
\begin{pf}
Let $\sigma\in\Lambda^{\Bbb N}$ be a sequence arbitrary. According to associative property of product of matrices, we can gain the following relation for each $X\in \Bbb {R}^{m}$;
 \begin{eqnarray*}
\begin{array}{ll}
 G_{\sigma_n}(X) & = AD_{\lambda_n}{\underbrace{A^{-1}\,A}_I}D_{\lambda_{n-1}}A^{-1}\ldots AD_{\lambda_2}{\underbrace{A^{-1}\,A}_I}D_{\lambda_1}A^{-1}\,X\\
                  & = AD_{\lambda_n}D_{\lambda_{n-1}}\ldots D_{\lambda_2}D_{\lambda_1}A^{-1}\,X
\end{array}
\end{eqnarray*}
and we have $F_{\sigma_n}(X)= D_{\lambda_n}D_{\lambda_{n-1}}\ldots D_{\lambda_2}D_{\lambda_1}\,X$.\\
Now, we define
$\left\{ \begin{array}{rl}
 h:\,\Bbb {R}^{m} & {\longrightarrow \Bbb {R}^{m}}\\
X & \longmapsto {AX.}
\end{array}\right.$
First notice that since $A$ is an invertible matrix so the function $h$ is a homeomorphism. For each $X\in \Bbb {R}^{m}$ we have:
\begin{eqnarray*}
\begin{array}{ll}
h\Big(F_{\sigma_n}(X)\Big) & = A\,F_{\sigma_n}(X)= A\,D_{\lambda_n}D_{\lambda_{n-1}}\ldots D_{\lambda_2}D_{\lambda_1}\,X\\
                           & = A\,F_{\sigma_n}(X)= A\,D_{\lambda_n}D_{\lambda_{n-1}}\ldots D_{\lambda_2}D_{\lambda_1}{\underbrace{A^{-1}\,A}_I}\,X\\
                           & =  G_{\sigma_n}\Big(A\,X\Big)= G_{\sigma_n}\Big(h(X)\Big).
\end{array}
\end{eqnarray*}
Thus $\mathcal{F}$ and $\mathcal{G}$ are weakly topological conjugates. $\Box$
\end{pf}
\section{Necessary Condition for Structural Stability IFSs}

  Now, we are going to define the concept structural stability for IFSs. In order to definition of distance of two IFS, we need to the following definitions that we state them from \cite{Pilyugin1999}.\\ Suppose $M$ is a $C^{\infty}$ smooth $m$-dimensional closed( that is, compact and boundaryless) manifold and $r$ is a Riemannian  metric on $M$. Let $f$ and $g$ be homeomorphisms on $M$; that is, $f, g\in Homeo(M)$. We define the metric $\rho_0$ as follows:
$$\rho_0(f, g)= Max\bigg\{r\Big(f(x), g(x)\Big),\,r\Big(f^{-1}(x),g^{-1}(x)\Big);\,\,\forall x\in M\bigg\}$$
Assume $f$ and $g$ are $C^{1}$-diffeomorphisms  on $M$; that is, $f, g\in Diff^{1}(M)$. We define the metric $\rho_1$ as follows;
$$\rho_1(f, g)= \rho_0(f, g)+ Max\bigg\{\parallel Df(x)-Dg(x)\parallel;\,\,\forall x\in M\bigg\};$$
 that here
 $${\scriptstyle Max\bigg\{\parallel Df(x)-Dg(x)\parallel;\,\,\forall x\in M\bigg\}= Max\bigg\{\mid Df(x)u-Dg(x)u\mid;\,\,\forall x\in M\;and\;\forall u\in T_{x}M:\, \mid u\mid= 1\bigg\}}$$
\begin{dfn}
Suppose $\mathcal{F}=\{f_{\lambda}, M : \lambda\in\Lambda\}$ and $\mathcal{G}=\{g_{\lambda},M : \lambda\in\Lambda\}$ are two IFSs as subsets of $Homeo(M)$. We denote the measure distance for two IFSs by ${\mathcal{D}}_0$ and define as follows:\\
If $\mathcal{F}=\mathcal{G}$ then put ${\mathcal{D}}_0\Big(\mathcal{F}, \mathcal{G}\Big)=0$\\
If $\mathcal{F}\neq\mathcal{G}$ then
$${\mathcal{D}}_0\Big(\mathcal{F}, \mathcal{G}\Big)= Max\Big\{\rho_0(f_{\lambda_i}, g_{\lambda_j}):\quad f_{\lambda_i}\in\mathcal{F}\; and\; g_{\lambda_j}\in\mathcal{G}\;for \;every\; \lambda\in\Lambda\;and\;i,j\in{\Bbb Z}\Big\}$$
If IFSs $\mathcal{F}$ and $\mathcal{G}$ are subset of $Diff^{1}(M)$ then we denote the measure distance for two IFSs by ${\mathcal{D}}_1$ and define as follows:\\
If $\mathcal{F}=\mathcal{G}$ then put ${\mathcal{D}}_1\Big(\mathcal{F}, \mathcal{G}\Big)=0$\\
If $\mathcal{F}\neq\mathcal{G}$ then
$${\mathcal{D}}_1\Big(\mathcal{F}, \mathcal{G}\Big)= Max\Big\{\rho_1(f_{\lambda_i}, g_{\lambda_j}):\quad f_{\lambda_i}\in\mathcal{F}\; and\; g_{\lambda_j}\in\mathcal{G}\;for \;every\; \lambda\in\Lambda\;and\;i,j\in{\Bbb Z}\Big\}$$
\end{dfn}
\begin{dfn}
Assume $\mathcal{F}=\{f_{\lambda}, M : \lambda\in\Lambda\}$ is an IFS where for every $\lambda\in\Lambda$,\,$f_{\lambda}\in Diff^{1}(M)$. We say that IFS\,$\mathcal{F}$ is $\emph{\textbf{structurally stable}}$ if for given $\epsilon>o$ there is $\delta>0$ such that for any IFS\,$\mathcal{G}=\{g_{\lambda}, M : \lambda\in\Lambda\}$
subset $Diff^{1}(M)$ with ${\mathcal{D}}_1\Big(\mathcal{F}, \mathcal{G}\Big)<\delta$ then for any $\sigma=\{\ldots,\lambda_{-1},\lambda_0,\lambda_1,\ldots\}$ and for every $n\in\Bbb Z$ there is a homeomorphism $h:\,M\rightarrow M$ with the following properties:\\
$\left \{\begin{array}{lll}
i)  & F_{\sigma_n}oh= hoG_{\sigma_n}, &  \\
ii) & r\Big(x, h(x)\Big)<\epsilon,   &  \forall x\in M.
\end{array}\right.$
\end{dfn}
Roughly speaking, an IFS is structurally stable if nearby IFSs are weakly topological conjugates with it; that is, nearby IFSs have qualitatively the same dynamics.\\
The next theorem states necessary conditions for structural stability IFSs.
\begin{thm}
If IFS $\mathcal{F}=\{f_{\lambda}, \Bbb R : \lambda\in\Lambda\}$ as subset of $Diff^{1}(\Bbb R)$ is structurally stable, then fixed points of the functions $\mathcal{F}$ are hyperbolic
\end{thm}
\begin{pf}
Assume that $\epsilon>0$ is given, then there exists $\delta>0$ by definition of structural stability IFSs. Consider IFS $\mathcal{G}=\{g_{\lambda},\Bbb R : \lambda\in\Lambda\}$ such that ${\mathcal{D}}_1\Big(\mathcal{F}, \mathcal{G}\Big)<\delta$. To get a contradiction suppose that there exists the function $f_{\lambda_i}, i\in\Bbb N,$ of $\mathcal{F}$ such that fixed point $p$ is not hyperbolic; that is, $\mid{f^{'}_{\lambda_i}}(p)\mid= 1$. Put $\sigma=\{\lambda_i, \lambda_i, \ldots\}$. According to structural stability of the IFS $\mathcal{F}$, for given $\sigma\in\Lambda^{\Bbb N}$ and $n=1$ there exists a homeomorphism $h$ on $\Bbb R$ such that $hof_{\lambda_i}=g_{\lambda_i}oh$. Clearly $\rho_1\Big(f_{\lambda_i},g_{\lambda_i}\Big)<\delta$ because $\rho_1\Big(f_{\lambda_i},g_{\lambda_i}\Big)<{\mathcal{D}}_1\Big(\mathcal{F}, \mathcal{G}\Big)<\delta$ by definition of the metric ${\mathcal{D}}_1$. As $\mathcal{G}$ is a arbitrary IFS so the function $g_{\lambda_i}\in Diff^{1}(\Bbb R)$ also is a arbitrary function, thus for given $\epsilon>0$ there is $\delta>0$ such that for every $g\in Diff^{1}(\Bbb R)$ that $\rho_1\Big(f_{\lambda_i},g\Big)<\delta$, there exists a homeomorphism $h$ on $\Bbb R$ that $hof_{\lambda_i}=goh$, this means that the function $f_{\lambda_i}$ is structurally stable. We know that if a diffeomorphism is structurally stable then its fixed points are hyperbolic, so the fixed points of the function $f_{\lambda_i}$ are hyperbolic and this is contradictory with $\mid{f^{'}_{\lambda_m}}(p)\mid= 1$ and the statement is proved. $\Box$
\end{pf}
\section{An outline of future challenges}
We may ask the questions  that have never been answered for example,\\
How we can define the limit sets and the limit points for an IFS?\\
What would happen if there wasn't the value of limit $\lim_{n\rightarrow +\infty}{{n_1}/{n_2}}$ or it was not zero, in the theorem\ref{mainTh2}?\\
Can we extend the concept of an IFS to the continuous systems and how can we generalize the Hartman-Grobman Theorem to these systems and how can we define the structural stability?

\section*{References}

\bibliography{mybibfile1}

\end{document}